\documentclass[12pt]{amsart}
\usepackage{amssymb}
\ifx\mathrm\undefined\let\mathrm\rm\fi
\ifx\mathbf\undefined\let\mathbf\bf\fi
\ifx\mathfrak\undefined\let\mathfrak\frak\fi
\ifx\mathcal\undefined\let\mathcal\cal\fi
\ifx\mathbb\undefined\let\mathbb\Bbb\fi
\ifx\emph\undefined\let\emph\it\fi
\newcommand{\g}{{{\mathfrak g}\,}}

\newcommand{\h}{{{\mathfrak h\,}}}

\newcommand{\Z}{{\mathbb Z}}

\newcommand{\R}{{\mathbb R}}
\newcommand{\C}{{\mathbb C}}

\newcommand{\Ref}[1]{{(\ref{#1})}}
\newcommand{\be}{\begin{displaymath}}
\newcommand{\ee}{\end{displaymath}}
\newcommand{\bea}{\begin{eqnarray*}}
\newcommand{\eea}{\end{eqnarray*}}

 \let\eps\varepsilon \let\epsilon\eps

\let\al\alpha

\let\la\lambda 
 
 \let\phi\varphi

\newcommand{\half}{\frac12}
\newcommand{\bean}{\begin{eqnarray}}
\newcommand{\eean}{\end{eqnarray}}

\newcommand{\vs}{\vspace{.5\baselineskip}}
\newtheorem%
{thm}{Theorem}%[section]
\newtheorem%
{proposition}[thm]{Proposition}
\newtheorem%
{lemma}[thm]{Lemma}
\newtheorem%
{lemmadef}[thm]{Lemma-Definition}
\newtheorem%
{corollary}[thm]{Corollary}
\newtheorem%
{conjecture}[thm]{Conjecture}
\title[Special Functions, Conformal Blocks, Bethe Ansatz, and $SL(3,\Z)$]
{Special Functions, Conformal Blocks, Bethe Ansatz, and $SL(3,\Z)$}
\author[G. Felder and A. Varchenko]
{G. Felder $^{\,\star}$ and
A.\ Varchenko$^{\,\diamond,1}$}
\thanks{$^1$ Supported in part by NSF grant DMS-9801582}

\begin{document}
\maketitle

\begin{center}
{\it
$^\star$ Departement Mathematik, ETH-Zentrum, 8092 Z\"urich, Switzerland,

felder@math.ethz.ch

\medskip

$^\diamond$Department of Mathematics, University of North Carolina
at Chapel Hill,

Chapel Hill, NC 27599-3250, USA,

av@math.unc.edu}
\end{center}

%\vsk1.5>
\centerline{November, 2000}
%\vsk1.8>

\begin{abstract}
This is the talk of the second author at the meeting "Topological Methods
in Physical Sciences", London, November 2000. We review our work on KZB equations.
\end{abstract}

%\vsk>
%\vsk0>

\thispagestyle{empty}

${}$ \\ \,% \newline   

{\bf KZB equations. }

%${}$ %\\ \, %\newline   

The Knizhnik-Zamolodchikov-Bernard (KZB)
equations are a system of differential equations
arising in conformal field theory on Riemann surfaces. 
For each $g, n\in \Z_{\geq 0}$, a simple complex Lie algebra $\g$, $n$ highest weight
$\g$-modules $V_i$ and a complex parameter $\kappa$, we have such a system of equations.
In the case of genus $g=1$, they have the form
\bean\label{KZB}
\kappa {\partial v\over \partial {z_j}}=-\sum_\nu h^{(j)}_\nu
{\partial v\over \partial{\la_\nu}}
+\sum_{l,\,l\neq j}r(z_j-z_l,\la,\tau)^{(j,l)}v,\qquad j=1,...,n,
\notag
\eean
\bean\label{heaT}
4\pi i \kappa {\partial v\over \partial\tau}\,=\,
\Delta_\la v+{1\over 2}\sum_{i,j}s(z,\la,\tau)^{(i,j)}v.
\notag
\eean
The unknown function $v$ takes values in the zero weight space $V[0]=
\cap_{x\in \h}\text{Ker}(x)$ of the tensor product $V=V_1\otimes \cdots \otimes V_n$
with respect to the Cartan subalgebra $\h$ of $\g$. It depends on variables $z_1,...,z_n\in \C$,
modulus $\tau$ of the elliptic curve 
and $\la=\sum \la_\nu h_\nu \in \h$, where $(h_\nu)$ 
is an orthonormal basis
of $\h$ with respect to a fixed invariant bilinear form. 
In the equation, $r,s\in \g\otimes \g$ are suitable given tensor valued functions,
\cite{FW}.

The second equation is called the KZB-heat equation.

{\bf Example.} For $\g=sl_N,\, n=1,\, V=S^{mN}\C^N,\, \h=\C^N/\C(1,...,1)$,
the weight-zero space $V[0]$ is one dimensional, 
the KZB equations are scalar equations $\partial v/\partial {z_1} =0$ and
\bean\label{CM-N}
4\pi i \kappa {\partial v\over \partial \tau }\,=\,
\sum_{i=1}^N{\partial^2 v\over \partial \lambda_i^2}
+2m(m+1)\sum_{1\leq i<j\leq N}\rho'(\la_i-\la_j,\tau)v\,.
\eean
Here
$'$ denotes the derivative with respect to the first argument and
$\rho$ is defined in terms of the first Jacobi theta function,
\bean
\theta(t,\tau)=-\sum_{j\in\Z}
e^{\pi i(j+\half)^2\tau+2\pi i(j+\half)(t+\half)}, \qquad
\rho(t,\tau)={\theta'(t,\tau)\over \theta(t,\tau)}.
\notag
\eean
Notice that $\rho'=-\wp+c$ where $\wp(t,\tau)$ is the  Weierstrass
function, $c=c( \tau)$ a function of $\tau$, and we recover in the right hand side 
of the KZB-heat 
equation the Hamilton operator $H_{N,m}$ of the elliptic Calogero-Moser quantum 
$N$-body system,
\bean
-H_{N,m}=\sum_{i=1}^N{\partial^2 \over \partial \lambda_i^2}
+2m(m+1)\sum_{ i<j}\rho'(\la_i-\la_j,\tau)\,
\notag
\eean
with coupling constant $m(m+1)$.

In particular, if $\g=sl_2$, then the  Cartan subalgebra
$\h$ can be identified
with $\C$, and $\la\in \C$. For the irreducible $2m+1$ 
dimensional module $V$, the KZB-heat equation takes the form
\bean\label{CM}
2\pi i \kappa {\partial v\over \partial \tau} 
(\la,\tau)\,=\,{\partial^2 v\over \partial \lambda^2}(\la,\tau)
\,+\, m(m+1) \rho'(\la,\tau)\,v(\la,\tau)\,.
\eean

${}$

A remarkable fact about all forms of the KZ equations is that they can be realized 
geometrically. They have solutions defined as hypergeometric integrals 
depending on parameters. 
% KZ equations can be realized as Gauss-Manin connections.

{\bf Example.} Consider the function
\bean
v(\la,\tau)\,=\,\int_0^1 \left({\theta (t,\tau)\over \theta '(0,\tau)}\right)
^{-{2\over\kappa}}
{\theta (\la-t,\tau) \theta '(0,\tau)\over \theta (\la,\tau) \theta (t,\tau)}
g(\la-{2\over \kappa}t,\tau)\,dt
\notag
\eean
where $g(\la,\tau)$ is any solution of the heat equation
\bean\label{}
2\pi i \kappa {\partial g\over \partial \tau}
(\la,\tau)\,=\,{\partial^2 g\over \partial \lambda^2}(\la,\tau).
\notag
\eean
For instance, $g(\la,\tau)=e^{\la\mu +{\mu^2\over 2\pi i \kappa}\tau}$ for  $\mu\in\C$.
Then $v$ is  a solution of the KZB-heat equation \Ref{CM} with $m=1$  \cite{FV1}.

{\bf Remark.} 
Simplest KZ equations have the classical Gauss hypergeometric function
as their solution. It is natural to consider solutions of all KZ type equations
as generalized hypergeometric functions. Thus there are hypergeometric functions
associated with curves of any genus. The function $v(\la,\tau)$ is the simplest 
elliptic hypergeometric function.

The Gauss hypergeometric function is
\bean
{\Gamma(b)\Gamma(c-b)\over \Gamma(c)}F(a,b,c;z)\,=\,
\int_1^\infty t^{a-c}(t-1)^{c-b-1}(t-z)^{-a}\,dt\,.
\notag
\eean
The functions 
\bean
\left({\theta (t,\tau)\over \theta '(0,\tau)}\right)^{a}\,\qquad
\text{and}\qquad
{\theta (\la-t,\tau) \theta '(0,\tau)\over \theta (\la,\tau) \theta (t,\tau)}
\notag
\eean
are elliptic analogs of the functions $t^a$ and $1/t$, respectively.
The hypergeometric solutions to the KZB equations were discovered through this analogy.

{\bf Remark.} If $\kappa$ tends to $0$, we are dealing with the eigenfunction problem.
Find eigenfunctions of the elliptic Calogero-Moser Hamiltonian,
\bean
\sum_{i=1}^N{\partial^2 v\over \partial \lambda_i^2}
+2m(m+1)\sum_{1\leq i<j\leq N}\rho'(\la_i-\la_j,\tau)\,v\,=\,E\,v\,.
\notag
\eean
The stationary phase method applied to hypergeometric solutions of the KZB equations 
gives eigenfunctions. For instance, application of the method to the function $v(\la,\tau)$
implies that for every $\mu\in \C$, the function $v(\la)=e^{\la\mu}\theta(\la-t_0,\tau)
/\theta(\la,\tau)$ is an eigenfunction of the operator
\bean
{\partial^2 \over \partial \lambda^2}
\,+\, 2 \rho'(\la,\tau)\,
\notag
\eean
if $t_0$ is a critical point of the function $e^{-\mu t}\theta(t,\tau)$.
This is a Bethe ansatz type formula, see  \cite{FV1, FV2}.
For $\g=sl_2$, the Bethe ansatz formulas reduce to Hermite's 1872 solution of the Lam\'e 
equation, see \cite{WW}.

${}$ %\\ \, \newline   

The fact that solutions have the form of explicitly written integrals is a  tool to study 
solutions as well as equations, for instance , the monodromy properties of
solutions or their modular properties with respect to changes of $\tau$.

${}$ %%%%%\\ \, \newline   

{\bf Modular symmetries of KZB.}

For any ${\scriptstyle\left(\begin{array}{cc}a&b\\ c&d\end{array}\right)}\in SL(2,\Z)$,
the elliptic curves with moduli $\tau$ and $(a\tau +b)/(c\tau + d)$ are isomorphic,
the corresponding KZB equations are related. Having a family of solutions one can ask
about monodromy properties of solutions with respect to transformations of the lattice
$\Z+\tau \Z$.

{\bf Example.} Consider equation \Ref{CM} for $m=0$,
\bean\label{=0}
2\pi i \kappa {\partial v\over \partial \tau} 
\,=\,{\partial^2 v\over \partial \lambda^2}\,.
\eean
If $v(\la,\tau)$ is a solution, then $\tilde v (\la,\tau) = v(\la,\tau + 1)$
and
\bean
\tilde v (\la,\tau) \,=\, {1\over \sqrt \tau} \, e^{-{\pi i\la^2\over 2\tau}}
\, v({\la\over \tau},\,-{1\over \tau})
\notag
\eean
are solutions too. If  $\{ v_\mu(\la,\tau)\}$ is  a family of solutions
depending on a parameter $\mu$, then one can ask about "monodromy" relations
of the three families of solutions: ${}$ $\{ v_\mu(\la,\tau)\}$,
$\{ v_\mu(\la,\tau +1)\}$, and $\{ {1\over \sqrt \tau} \, e^{-{\pi i\la^2\over 2\tau}}\,
v_\mu({\la\over \tau},\,-{1\over \tau})\}$.

{\bf The KZB-heat equation as a flat connection.} One can consider equation \Ref{CM} as 
the equation for horizontal sections of a connection on a vector bundle
over the upper half plane
$\frak H_+$ whose fiber $F(\tau)$ over a point $\tau$ is the space of functions
of $\la$.

%{\bf Sub-bundle of conformal blocks.} 
If $\kappa$ is a positive integer not less than $2$,
then the bundle has a finite dimensional sub-bundle (of
conformal blocks) invariant with respect to the connection
and consisting of certain theta functions of level $\kappa$, \cite{EFK, FW, FV1}.

${}$

In this paper we address the following 

${}$

{\bf Problem.} {\it Quantize the KZB heat equation, study modular properties of
the quantization.}

It turns out that the KZB-heat equation is quantized to a difference
qKZB-heat equation (with step $p$) in such a way that the $SL(2,\Z)$ symmetry
of the KZB-heat equation 
related to the lattice $Z+\tau\Z$ is quantized to an $SL(3,\Z)$ symmetry
related to the lattice $Z+\tau\Z+ p\Z$.

${}$

${}$

A quantization of the heat equation is a discrete connection over $\frak {H}_+$
with the same fiber, i.e. a linear operator $T(\tau,\tau +p): F(\tau + p)\to
F(\tau )$. This linear operator  tends to 
\bean
1\,+\,\text{const}\, p\, \, (\,{\partial^2 \over \partial \lambda^2}
\,+\, m(m+1) \rho'(\la,\tau) \, )\,+\,...
\notag
\eean
as $p\to 0$.

${}$

%{\bf Example.} 
For methodological reasons we describe first a quantization 
of  equation \Ref{=0}. For a nonzero $\eta \in \C$, introduce linear
operators
\bean
U\,&:&\, f(\la)\,\mapsto\, \,{i\over \sqrt {4 i\eta}}
\int_{\eta \R} e^{-\pi i {\la \mu \over 2 \eta}}\,f(-\mu)\,d\mu\,,
\notag
\\
\alpha\,&:&\, f(\la)\,\mapsto\,e^{-\pi i {\la^2\over 4 \eta}}\,f(\la)\,.
\notag
\eean
Define
\bean\label{im=0}
T(\tau,\tau + p)=\al U \al\,:\,
f(\la)\,\mapsto\, \,{i\over \sqrt {4 i\eta}}\,
\int_{\eta \R} e^{-\pi i {(\la +\mu)^2 \over 4 \eta}}\,f(-\mu)\,d\mu\,.
\eean
The translation operator $T(\tau,\tau + p)$ in this case is basically 
just the Fourier transform.

The qKZB-heat equation is the equation for flat sections,
\bean\label{heat=0}
v(\la,\tau)\,=\, \,{i\over \sqrt {4 i\eta}}\,
\int_{\eta \R} e^{-\pi i {(\la +\mu)^2 \over 4 \eta}}\,
v(-\mu,\tau + p)\,d\mu\,.
\eean
Set $p=-2\kappa\eta$. Let $\eta \to 0$ and $\la, \tau, \kappa$  fixed.

${}$

{\bf Theorem.} {\it Let $v_\eta(\la,\tau)$ be a family of solutions of \Ref{heat=0}
with asymptotics $v_\eta(\la,\tau)=v^0(\la,\tau)+ \eta v^1(\la,\tau)+...$ 
Then $v^0(\la,\tau)$ satisfies \Ref{=0}.}

${}$

{\bf Proof.} The stationary phase asymptotic expansion of the right hand side
of \Ref{heat=0} is 
\bean
v^0(\la,\tau)+ \eta v^1(\la,\tau)+ {i\eta\over \pi}
(2\pi i \kappa {\partial v^0\over \partial \tau}(\la,\tau) 
- {\partial^2 v^0\over \partial \lambda^2}(\la,\tau)) + O(\eta^2)+... 
\notag
\eean

$\square$

Introduce a function 
\bean\label{um=0}
u(\la,\mu,\tau,p,\eta)=e^{-\pi i {\la\mu\over 2\eta}}\,.
\eean
For every $\mu \in\C$, the function $u$ is a projective solution
of the qKZB-heat equation \Ref{heat=0},
\bean
u(\la,\mu,\tau,p,\eta)=e^{-\pi i {\mu^2\over 4\eta}}\,(T(\tau,\tau + p) u)
(\la,\mu,\tau +p,p,\eta),
\notag
\eean
in particular, ${v}(\la,\mu,\tau,p,\eta)=e^{\pi i {\mu^2\over 4\eta}{\tau\over p}
-\pi i {\la\mu\over 2\eta}}$ is a true solution.

${}$

{\bf Remark} \cite{FV3}.  Let $\kappa$ be a positive integer.
The  functions
\bean
\theta_{j,\kappa}(\lambda,\tau)=\sum_{r\in\Z+j/2\kappa}
e^{2\pi i\kappa (r^2\tau+r\lambda)},\qquad j\in\Z/2\kappa\Z,
\notag
\eean
form the space   $\Theta_{\kappa}(\tau)$  of theta functions
of level $\kappa$. Let $E(\tau)=\{f\in \Theta_{\kappa}(\tau)\,|\,
f(-\la)=-f(\la)\}$ be the space of odd theta functions. 
%The space $E(\tau)$ is a finite dimensional subspace of the fiber $F(\tau)$.
The translation operator $T(\tau,\tau+p)$ maps $E(\tau+p)$ to
$E(\tau)$ if $-p/2\eta=\kappa$.

 This statement is based on the identity
\bean
\theta_{j,\kappa}(\lambda,\tau)=\frac i{\sqrt{4i\eta}}
\int_{2\eta\R}e^{-\frac{i\pi}{4\eta}(\lambda+\mu)^2}
\theta_{j,\kappa}(-\mu,\tau-2\eta\kappa)\,d\mu\,.
\notag
\eean
The space $E(\tau)$ is the quantization of the finite dimensional space of conformal blocks.

${}$

${}$

{\bf Modular properties of the q-heat operator (in this case of the Fourier transform).}

 The group $SL(3,\Z)$  is generated
by the elementary matrices $e_{i,j}$, $i\neq j$.
The elementary matrix $e_{i,j}$ is
the element of $SL(3,\Z)$ which differ from the identity
matrix by having the $i,j$ matrix element equal to $1$.
 The relations can be chosen  \cite{M} to be
\begin{eqnarray*}
e_{i,j}e_{k,l}&=&e_{k,l}e_{i,j},\qquad i\neq l,\quad j\neq k,\\
e_{i,j}e_{j,k}&=&e_{i,k}e_{jk}e_{i,j},\\
(e_{1,3}e_{3,1}^{-1}e_{1,3})^4&=&1.
\end{eqnarray*}
Consider $\C^3$ with coordinates $x_1,x_2,x_3$. 
The group $SL(3,\Z)$ acts on $\C^3$ in the  standard way.
Consider the trivial bundle over $\C^3$ with the same fiber $F$ over a point $x\in\C^3$.
We define a projectively flat connection over the orbit of a point.
Namely, for any generator $e_{i,j}$, we define a linear operator
\bean
\phi_{i,j}(x)\,:\, F(e^{-1}_{i,j}x)\,\to\, F(x)
\notag
\eean
so that all relations in $SL(3,\Z)$ are projectively satisfied. That means that
the linear operator corresponding to the left hand side of a relation is equal
to the linear operator corresponding to the right hand side of the relation multiplied
by a number.

${}$

{\bf Remark.} The operator $\phi_{2,1}(p,\tau,\eta)\,:\, F(p,\tau - p, \eta)\,\to\, 
F(p,\tau, \eta)$ will correspond to the qKZB-heat operator.

${}$

Introduce linear operators
\bean
U(x_1,x_2,x_3)\,&:&\, f(\la)\,\mapsto\, 
%-\,{1\over 4\pi \sqrt {i x_3}}
\int_{ x_3\R} e^{-\pi i {\la \mu \over 2  x_3}}\,f(-\mu)\,d\mu\,,
\notag
\\
\alpha(x_3)\,&:&\, f(\la)\,\mapsto\,e^{-\pi i {\la^2\over 4 x_3}}\,f(\la)\,,
\notag
\\
\beta(x_1,x_2,x_3)\,&:&\, f(\la)\,\mapsto\,e^{-\pi i {\la^2\over 4 }{x_1\over x_2x_3}}
\,f(\la)\,.
\notag
\eean
Set
\bean
\phi_{1,3}(x_1,x_2,x_3)&=&1,
\notag
\\
 \phi_{2,3}(x_1,x_2,x_3)&=& 1,
\notag
\\
\phi_{1,2}(x_1,x_2,x_3)&=&\alpha (x_3), 
\notag
\\
\phi_{3,2}(x_1,x_2,x_3)&=&\beta (x_1,x_2-x_3,x_3),
\notag
\\
\phi_{2,1}(x_1,x_2,x_3)&=&\alpha (x_3)\,U(x_2,x_2-x_1,x_3)\,\alpha (x_3),
\notag
\\
\phi_{3,1}(x_1,x_2,x_3)&=&\beta (x_1-x_3,-x_3,x_2)^{-1}\,
U(x_1-x_3,-x_3,x_2)^{-1}\,\beta (x_3,x_2,x_3-x_1)^{-1}.
\notag
\eean

${}$

{\bf Theorem.} {\it The operators $\phi_{i,j}$ define a projectively flat connection
over orbits of the $SL(3,\Z)$ action.}

${}$

 {\bf Remark.} Consider $\C^3$ and the projectivization of the dual space,
$P(\C^3)^*$. Consider $X\subset (\C-0)^3\times P(\C^3)^*$ where
\bean
X=\{((x_1,x_2,x_3),(y_1:y_2:y_3))\,|\, x_1y_1+x_2y_2+x_3y_3=0\,\}\,.
\notag
\eean
The natural projection $X\to \C^3-0$ is a projective line bundle. The group $SL(3,\Z)$
acts on $X$, $g:(x,y)\mapsto (gx, (g^t)^{-1}y)$ for any $g\in SL(3,\Z)$. Fix an affine coordinate 
on fibers, $t:(y_1:y_2:y_3)\mapsto y_2/y_1$. Then for any $i,j$, $i\neq j$,
we have
\bean
e_{i,j}:(e^{-1}_{i,j} x,t)\mapsto (x,f_{i,j}(x,t))
\notag
\eean
where $f_{1,3}(x,t)=f_{2,3}(x,t)=t$, $f_{1,2}(x,t)=t-1$,
$f_{3,2}(x,t)= (tx_3+x_1)/(x_3-x_2)$, $f_{2,1}(x,t)=t/(1-t)$,
$f_{3,1}(x,t)=t(x_3-x_1)/(tx_2+x_3)$.

The $SL(3,\Z)$-action on $X$ is closely related to the projectively flat connection
described in the Theorem. Namely, set
\bean
G(\lambda;x_1,x_2,x_3;t) \,= \,e^{{\pi i\over 4 }\lambda^2{t\over x_3}}\,.
\notag
\eean
Call a Gaussian in the fiber $F(x_1,x_2,x_3)$  a function of the form
const$\cdot G(\lambda;x_1,x_2,x_3;t) $ for some number $t$. The linear operators 
$ \phi_{i,j}(x)\,:\, F(e^{-1}_{i,j}x)\,\to\, F(x)$ preserve the Gaussians.
Moreover, for all $i,j$, we have
\bean
\phi_{i,j}(x)\,:\,G(\lambda;x_1,x_2,x_3;t)
\mapsto \text{const}\cdot G(\lambda;x_1,x_2,x_3;f_{i,j}(t,x))\,.
\notag
\eean

${}$

 {\bf Remark.} The Theorem easily follows from the following two main equations satisfied by the
Fourier transform, which we call

the {\it  q-heat equation,}
\bean\label{q-heat}
\alpha(x_3)\,U(x_1,x_1+x_2,x_3)\,\alpha(x_3)
\,U(x_1+x_2,x_2,x_3)\,\alpha(x_3)\,=\,U(x_1,x_2,x_3)\,,
\eean

and the {\it modular equation,}
\bean\label{modular}
U(x_3,x_2,-x_1)\,\beta (-x_3,x_2,x_1)\,U(x_1,-x_3,x_2)\,=
\\
\,\beta(x_2,x_3,x_1)\,U(x_1,x_2,x_3)\,\beta(x_1,x_2,x_3)\,.
\notag
\eean

${}$

{\bf Proof of the q-heat equation.}
\bean
\,e^{-\pi i {\la^2\over 4 x_3}}\,
%{1\over 4\pi \sqrt {i x_3}}
\int  e^{-\pi i {\la \nu \over 2  x_3}}
\left(\,e^{-\pi i {\nu^2\over 4 x_3}}
%{1\over 4\pi \sqrt {i x_3}}
\int  e^{\pi i {\nu \mu \over 2  x_3}}\,e^{-\pi i {\mu^2\over 4 x_3}}\,f(-\mu)\,d\mu\,
\right)\,d\nu\,=
\notag
\\
\int \,\left(\, \int e^{-\pi i {(\la +\mu-\nu)^2 \over 4  x_3}} d\nu\,\right)\,
e^{-\pi i {\la \mu \over 2  x_3}}\,f(-\mu)\,d\mu\,=
\notag
\\
\text{const}(x_3)\,
%-\,{1\over 4\pi \sqrt {i x_3}}
\int e^{-\pi i {\la \mu \over 2  x_3}}\,f(-\mu)\,d\mu\,.
\notag
\eean
{\bf Proof of the modular equation.}
\bean
%{1\over 4\pi \sqrt {-i x_1}}
\int  e^{\pi i {\la \nu \over 2  x_1}}
e^{\pi i {\nu^2\over 4 }{x_3\over x_1x_2}}\,
%{1\over 4\pi \sqrt {i x_2}}
\left(\,\int  e^{\pi i {\nu \mu \over 2  x_2}}\,f(-\mu)\,d\mu\,
\right)\,d\nu\,=
\notag
\\
\int \,\left(\,\int \,
e^{\pi i {1\over 4 x_1x_2x_3}( x_1\mu + x_2\la  +  x_3\nu)^2}
\, d\nu\,\right)\,
e^{-{\pi i\over 4}( \la^2 {x_2\over x_1x_3}
+2{\la\mu\over x_3}+
\mu^2{x_1\over x_2x_3})}
f(-\mu)\,d\mu\, =
\notag
\\
\text{const}(x_1,x_2,x_3)\,
%-\,{1\over 4\pi \sqrt {i x_3}}
e^{-\pi i {\la^2\over 4 }{x_2\over x_1x_3}}
\int e^{-\pi i {\la \mu \over 2  x_3}}\,
e^{-\pi i {\mu^2\over 4 }{x_1\over x_2x_3}}\,
f(-\mu)\,d\mu\,.
\notag
\eean

${}$

${}$

Similarly, one can quantize equation \Ref{CM} for $m=1$,
\bean\label{m=1}
2\pi i \kappa {\partial v\over \partial \tau} 
(\la,\tau)\,=\,{\partial^2 v\over \partial \lambda^2}(\la,\tau)
\,+\, 2 \rho'(\la,\tau)\,v(\la,\tau)\,.
\eean

Introduce a function
\bean\label{um=1}
u(\la,\mu,\tau,p,\eta)\,=\,e^{-\pi i {\la\mu\over 2\eta}}\,\int_0^1 \,
\Omega_{2\eta}(t,\tau,p)\,
{\theta(\la+t,\tau) \theta (\mu+t,p)\over \theta(t-2\eta,\tau) \theta (t-2\eta,p)}\,
dt
\eean
where
\bean
\Omega_a(t,\tau,p)=\prod_{j,k=0}^\infty \frac {(1-e^{2\pi
    i(t-a+j\tau+kp)})(1-e^{2\pi i(-t-a+(j+1)\tau+(k+1)p)})}
{(1-e^{2\pi i(t+a+j\tau+kp)})(1-e^{2\pi i(-t+a+(j+1)\tau+(k+1)p)})}\,.
\notag
\eean
This is the analog for $m=1$ of  function \Ref{um=0}.

Define the translation operator for $m=1$ as
\bean
T(\tau,\tau + p)\,:\, f(\la)\,\mapsto \,
-\,{1\over 4\pi \sqrt {i \eta}}\, e^{-\pi i {\la^2\over 4\eta}}\,
\int_{ \eta\R} 
u(\la,\mu,\tau,\tau+p,\eta)\,&\times
\notag
\\
{\theta(4\eta,\tau+p) \theta' (0,\tau + p)\over \theta(\mu-2\eta,\tau+p) 
\theta (\mu + 2\eta, \tau +p)}\,
&e^{-\pi i {\mu^2\over 4\eta}}\,f(-\mu)\, d\mu\,.
\notag
\eean
%This is the analog for $m=1$ of the operator \Ref{im=0}.

The qKZB-heat equation is the equation
\bean\label{qm=1}
v(\la,\tau)\,=\,(T(\tau,\tau + p)v)(\la,\tau + p )\,.
\eean

${}$

It turns out that the translation operator for $m=1$ has properties analogous
to the properties of the Fourier transform, see \cite{FV3, FV4, FV5}.
Namely,

{\it 
\begin{enumerate}

\item The semiclassical limit of equation \Ref{qm=1} is equation \Ref{m=1}.

\item  For every $\mu \in\C$, the function $u$ in \Ref{um=1}
is a projective solution of the qKZB-heat equation \Ref{qm=1},
\bean
u(\la,\mu,\tau,p,\eta)\,=\,e^{-\pi i {\mu^2\over 4\eta}}\,(T(\tau,\tau + p) u)
(\la,\mu,\tau +p,p,\eta),
\notag
\eean
in particular, ${v}(\la,\mu,\tau,p,\eta)=e^{\pi i {\mu^2\over 4\eta}{\tau\over p}}
u(\la,\mu,\tau,p,\eta)$ is a true solution.

\item The integral operator $T$ has an $SL(3,\Z)$ symmetry similar to the $SL(3,\Z)$ 
symmetry of the Fourier transform.
\end{enumerate}
}

${}$

{\bf Remark.} The fact that $u$ is a ( projective ) solution of the qKZB-heat equation
is a relation  of the form
\bean\label{uu}
u(\la,\nu,\tau,p)\,=\,u(\la,\mu,\tau,\tau + p)*u(\mu,\nu,\tau + p,p)
\eean
where $*$ is a suitable convolution. If $\tau, p$ tend to infinity,
then the function $u$ has a trigonometric limit.
In this limit, equation \Ref{uu} becomes a simplest example of the
Macdonald-Mekhta identity which has the form
$u(\la,\nu)\,=\,u(\la,\mu)*u(\mu,\nu)$.
The theory of the
qKZB-heat equation is an elliptic analog of Macdonald's theory, see \cite{EV1}. 

${}$

{\bf Remark.} We have
\bean
\Omega_a(t,\tau,p)\,=\,{\Gamma (t+a,\tau,p)\over \Gamma (t-a,\tau,p)}
\notag
\eean
where
\bean
\Gamma(t,\tau,p)=\prod_{j,k=0}^\infty  { 1-e^{2\pi i(-t+(j+1)\tau+(k+1)p)} 
\over
1-e^{2\pi i(t+j\tau+kp)} }\,
\notag
\eean
is the elliptic gamma function. The elliptic gamma function
has an $SL(3,\Z)$ symmetry \cite{FV5} based on  the following two main equations, which we call

the {\it  q-heat equation,}
\bean\label{gamma-heat}
\Gamma(t+\tau,\tau,\tau + p)\,{}\,\Gamma(t,\tau+p,p)\,{}\,=\,{}\,
\Gamma(t,\tau,p),
\notag
\eean

and the {\it modular equation,}
\bean\label{gamma-modular}
\Gamma( {t\over p},{\tau\over p},-{1\over p})\,{}\,=\,{}\,
e^{i\pi Q(t;\tau,p)}\,{}\,
\Gamma({t-p\over \tau},-{1\over \tau},-{p\over\tau})\,{}\,
\Gamma(t,\tau,p),
\notag
\eean
where
\begin{eqnarray*}
Q(t;\tau,p)&=&
\frac{1}{3\tau p}\,t^3
-
\frac{\tau+p-1}{2\tau p}
t^2
+
\frac{\tau^2+p^2+3\tau p-3\tau-3p+1}
{6\tau p}
t
\\ & &
+
\frac1{12}
(\tau+p-1)
(\tau^{-1}+p^{-1}-1),
\end{eqnarray*}
cf. \Ref{q-heat}, \Ref{modular}.

The q-heat and modular equations for the gamma function imply the $SL(3,\Z)$
symmetry of the translation operator for $m=1$.

${}$

${}$

${}$

\end{document}